\documentclass[titlepage,12pt]{article}
\usepackage{amssymb}
\usepackage{amsfonts}
\usepackage{amsmath}
\begin{document}

{\Large \bf An Ergodic Result} \\ \\

{\bf Elem\'{e}r E Rosinger} \\
Department of Mathematics \\
and Applied Mathematics \\
University of Pretoria \\
Pretoria \\
0002 South Africa \\
eerosinger@hotmail.com \\ \\

{\bf Abstract} \\

A rather general ergodic type scheme is presented on arbitrary sets $X$, as they are generated
by arbitrary mappings $T : X \longrightarrow X$. The structures considered on $X$ are given by
suitable subsets of the set of all of its finite partitions. Ergodicity is studied not with
respect to subsets of $X$, but with the {\it inverse limits} of families of finite
partitions. \\ \\

{\bf 1. The Setup} \\

Let $( X, \Sigma, T)$ be as follows : $X$ is an arbitrary nonvoid set, $\Sigma \subseteq
{\cal P} ( X )$ is a nonvoid set of subsets of $X$, while $T : X \longrightarrow X$. \\
The issue considered, as usual in Ergodic Theory, is the behaviour of the sequence of iterates
$T^n ( x ),~ n \in \mathbb{N}_+ = \{ 1, 2, 3, \ldots \}$, for an arbitrary given $x \in X$. Of
a main interest in this regard is of course the case when $X$ is infinite. \\

A simplest and natural way to follow is to consider a {\it partition} of $X$, and see how the
mentioned sequence of iterates may possibly move through the various sets of that partition.
In this regard, a further simplest and natural case is when the partitions considered for $X$
are {\it finite}, and thus at least one of their sets must contain {\it infinitely} many terms
of any such sequence of iterates. \\

As it turns out, a number of properties can be obtained simply form the finite versus infinite
{\it interplay} as set up above, an interplay slightly extending the usual pigeon-hole
principle. However, in order to obtain such properties, one may have to {\it shift} the usual
focus which tends to be concerned with the relationship between the mapping $T$ and its
iterates $T^n$, with $n \in \mathbb{N}_+$, and on the other hand, the various subsets $A
\subseteq X$. Namely, this time one is dealing with the relationship between the mapping $T$
and its iterates $T^n$, with $n \in \mathbb{N}_+$, and on the other hand, whole {\it families}
of {\it finite} partitions $\Delta$ of $X$. \\

Let us therefore consider \\

(1.1)~~~ $ {\cal FP} ( X, \Sigma ) $ \\

the set of all {\it finite} partitions $\Delta$ of $X$ with nonvoid subsets in $\Sigma$, thus
$\Delta \subseteq \Sigma$, $\Delta$ is finite, and $X = \bigcup_{A \,\in\, \Delta} A$, where
for $A \in \Delta$ we have $A \neq \phi$, however in general, none of $A \in \Delta$ need to
be finite. \\

Given $x \in X$ and $\Delta \in {\cal FP} ( X, \Sigma )$, then obviously \\

(1.2)~~~ $ \exists~~ A \in \Delta ~:~
                  \{~ n \in \mathbb{N}_+ ~|~ T^n ( x ) \in A ~\} ~~\mbox{is infinite} $ \\

since $\Delta$ is finite. \\

Let us therefore denote \\

(1.3)~~~ $ \Delta( x ) = \{~ A \in \Delta ~|~
                  \{~ n \in \mathbb{N}_+ ~|~ T^n ( x ) \in A ~\} ~~\mbox{is infinite} ~\} $ \\

and then (1.2) implies \\

(1.4)~~~ $ \Delta ( x ) \neq \phi $ \\

{\bf Problem 1} \\

Given $x \in X$, what happens with $\Delta( x )$, when $\Delta$ ranges over \\ ${\cal FP} ( X,
\Sigma )$ ? \\

{\bf Example 1} \\

Let $X = \mathbb{N},~ \Sigma = {\cal P} ( \mathbb{N} )$ and consider the following three
cases of mappings $T : \mathbb{N} \longrightarrow \mathbb{N}$, where here and in the sequel,
we denote $\mathbb{N} = \{ 0, 1, 2, 3, \ldots \}$ : \\

1) $T$ is given by the usual shift $T ( x ) = x + 1,~ x \in \mathbb{N}$. \\

If $\Delta \in {\cal FP} ( \mathbb{N}, {\cal P} ( \mathbb{N} ) )$, then obviously there exists
$A \in \Delta$ such that $A$ is infinite. Furthermore, for every $x \in \mathbb{N}$, we
have \\

(1.5)~~~ $ \Delta ( x ) = \{~ A \in \Delta ~|~ A ~~\mbox{is infinite} ~\} \neq \phi $ \\

2) If $T$ is the identity mapping then clearly \\

(1.6)~~~ $ \Delta ( x ) = \{~ A ~\} \neq \phi,~~~\mbox{where}~ x \in A \in \Delta $ \\

3) Let us now assume that, for a given $x_* \in \mathbb{N}$, we have $T ( x ) = x_*$, with
$x \in \mathbb{N}$. Then obviously \\

(1.7)~~~ $ \Delta ( x ) = \{~ A ~\} \neq \phi,~~~\mbox{where}~ x_* \in A \in \Delta $ \\

{\bf Remark 1} \\

The above general setup clearly contains as a particular case the following one which is of a
wide interest in Ergodic Theory, namely, $( X, \Sigma, \nu )$, where $\Sigma$ is a
$\sigma$-algebra on $X$, while $\nu$ is a probability on $( X, \Sigma )$. In that case, the
mapping $T$ is supposed to satisfy the conditions \\

(1.8)~~~ $ T^{-1} ( \Sigma ) \subseteq \Sigma $ \\

and \\

(1.9)~~~ $ \nu ( T^{-1} ( A ) ) = \nu ( A ),~~~ A \in \Sigma $ \\

We note that $\Sigma$ being a $\sigma$-algebra, we have in particular \\

(1.10)~~~ $ \begin{array}{l}
                    ~~~~*)~~ \forall~~ A, A\,' \in \Sigma ~:~ A\,'' = A \cap A\,' \in \Sigma \\ \\
                    ~~~**)~~ \phi,~ X \in \Sigma
             \end{array} $ \\

consequently \\

(1.11)~~~ $ {\cal FP} ( X, \Sigma ) \neq \phi $ \\ \\

{\bf 2. Towards a Solution} \\

First we observe the following natural structure on ${\cal FP} ( X, \Sigma )$, given by the
concept of {\it refinement}. Namely, if $\Delta , \Delta\,' \in {\cal FP} ( X, \Sigma )$, we
define \\

(2.1)~~~ $ \Delta \leq \Delta\,' $ \\

if and only if \\

(2.2)~~~ $ \forall~~ A\,' \in \Delta\,' ~:~
                        \exists~~ A \in \Delta ~:~ A\,' \subseteq A $ \\

and in view of that, we can define the mapping \\

(2.3)~~~ $ \psi_{\Delta\,',\, \Delta} : \Delta\,' \longrightarrow \Delta\ $ \\

by \\

(2.4)~~~ $ A\,' \subseteq A = \psi_{\Delta\,',\, \Delta} ( A\,' ),~~~ A\,' \in \Delta\,' $ \\

Then we obtain \\

{\bf Lemma 1} \\

(2.5)~~~ $ \psi_{\Delta\,',\, \Delta} ( \Delta\,'( x ) )
                            \subseteq \Delta( x ),~~~ x \in X $ \\

{\bf Proof} \\

If $A\,' \in \Delta\,'( x )$, then (1.3) gives \\

$~~~~~~ \{~ n \in \mathbb{N}_+ ~|~ T^n ( x ) \in A\,' ~\} ~~\mbox{is infinite} $ \\

but in view of (2.4), we have \\

$~~~~~~ A\,' \subseteq \psi_{\Delta\,',\, \Delta} ( A\,' ) $ \\

hence \\

$~~~~~~ \{~ n \in \mathbb{N}_+ ~|~ T^n ( x ) \in \psi_{\Delta\,',\, \Delta} ( A\,' ) ~\}
                                                                  ~~\mbox{is infinite} $ \\

thus (2.5).

\hfill $\Box$ \\

Let us pursue the consequences of the above result in (2.5). In this regard we note that in
the usual particular case in Remark 1, the partial order (2.1) on ${\cal FP} ( X, \Sigma )$
is in fact {\it directed}, and obviously has the following stronger property \\

(2.6)~~~ $ \forall~~ \Delta,~ \Delta\,' \in {\cal FP} ( X, \Sigma ) ~:~
              \exists~~ \Delta \vee \Delta\,' \in {\cal FP} ( X, \Sigma )$ \\

since \\

(2.7)~~~ $ \Delta \vee \Delta\,' = \{~ A \cap A\,' ~|~
               A \in \Delta,~ A\,' \in \Delta\,',~ A \cap A\,' \neq \phi ~\} $ \\

However, for a greater generality, let us consider in Problem 1 not only the whole of ${\cal
FP} ( X, \Sigma )$, but also arbitrary subsets of it. Let therefore $( \Lambda, \leq )$ be any
partially ordered set, and consider a mapping \\

(2.8)~~~ $ \Lambda \ni \lambda \longmapsto \Delta_\lambda \in {\cal FP} ( X, \Sigma ) $ \\

such that \\

(2.9)~~~ $ \lambda \leq \lambda\,' \Longrightarrow
                          \Delta_\lambda \leq \Delta_{\lambda\,'} $ \\

We call the family $( \Delta_\lambda )_{\lambda \in \Lambda}$ a {\it refinement chain}. \\

Obviously, in view of the above, ${\cal FP} ( X, \Sigma )$ itself is such a refinement chain,
namely, with $\Lambda = {\cal FP} ( X, \Sigma )$, the partial order in (2.1), and with the
identity mapping in (2.8). \\

The main point to note is the following. Given $x \in X$, then (1.4) implies \\

(2.10)~~~ $ \Delta_\lambda ( x ) \neq \phi,~~~ \lambda \in \Lambda $ \\

Hence in view of (2.5), (2.9), we have for $\lambda \leq \lambda\,'$ \\

(2.11)~~~ $ \phi \neq \psi_{\Delta_{\lambda\,'},\, \Delta_\lambda }
                          ( \Delta_{\lambda\,'} ( x ) )\subseteq \Delta_\lambda ( x ) $ \\

Now, based on (2.10), let us use the notation \\

(2.12)~~~ $ \Delta_\lambda ( x ) =
                  \{~ A_{\lambda,1}(x), \ldots , A_{\lambda,m_\lambda}(x) ~\} $ \\

where $m_\lambda \geq 1$, and $\phi \neq A_{\lambda,j}(x) \in \Sigma$, with $1 \leq j
\leq m_\lambda$. \\

{\bf Problem 2} \\

A more precise reformulation of Problem 1 is as follows. We can investigate whether for a
given $x \in X$, one or the other of the following two properties may hold, namely \\

(2.13)~~~ $ \exists~~ \Lambda \ni \lambda \longmapsto A_\lambda \in \Delta_\lambda(x)  ~:~
                                             \bigcap_{\, \lambda \in \Lambda}\, A_\lambda \neq \phi $ \\

or what appears to be a milder property \\

(2.14)~~~ $ \underleftarrow{\lim}_{\, \lambda \,\in\, \Lambda}~ \Delta_\lambda(x) \neq \phi $

\hfill $\Box$ \\

The {\it inverse limit}, [1, p. 191], in (2.14) is of the family \\

(2.15)~~~ $ ( \Delta_\lambda(x) ~|~ \lambda \in \Lambda ) $ \\

with the mappings, see (2.8), (2.10) \\

(2.16)~~~ $ \psi_{\lambda\,',\, \lambda,\, x} : \Delta_{\lambda\,'}( x )
                                          \longrightarrow \Delta_\lambda( x ) $ \\

for $\lambda, \lambda\,' \in \Lambda,~ \lambda \leq \lambda\,'$ where \\

(2.17)~~~ $ \psi_{\lambda\,',\, \lambda} : \Delta_{\lambda\,'}
                                            \longrightarrow \Delta_\lambda $ \\

is given by \\

(2.18)~~~ $ \psi_{\lambda\,',\, \lambda} = \psi_{\Delta_{\lambda\,'},\, \Delta_\lambda } $ \\

while \\

(2.19)~~~ $ \psi_{\lambda\,',\, \lambda,\, x} =
                \psi_{\lambda\,',\, \lambda} \,|\,_{{}_{\Delta_{\lambda\,'}}(x)} $ \\

In order to establish (2.14), we recall the definition of the inverse limit, namely \\

(2.20)~~~ $ \begin{array}{l}
              \underleftarrow{\lim}_{\, \lambda \,\in\, \Lambda}~ \Delta_\lambda(x) = \\ \\
                ~~~= \left \{~ ( A_\lambda ~|~ \lambda \in \Lambda ) \in
                                  \prod_{\,\lambda \,\in\, \Lambda} \Delta_\lambda(x) ~~~
                       \begin{array}{|l}
                            ~\forall~~ \lambda, \lambda\,' \in
                                         \Lambda, \lambda \leq \lambda\,' ~: \\
                             ~~~~~ \psi_{\lambda\,',\, \lambda,\, x}(A_{\lambda\,'}) =
                                                  A_\lambda
                       \end{array}  ~\right \}
            \end{array} $ \\

We note that in the definition of the inverse limit, the partial order $\leq$ on $\Lambda$ can
be arbitrary, and in fact, it can be a mere pre-order. \\

Further we note, [1, Exercise 4, no. 4, p. 252], that an inverse limit such as for instance
in (2.20), can be void even when all sets $\Delta_\lambda(x)$ are nonvoid and all mappings
$\psi_{\lambda\,',\, \lambda,\, x}$ are surjective. \\

However, as seen in Theorem 1 in the sequel, this is not the case in (2.14). \\

Meanwhile, for the sake of further clarification, we consider (2.20) in the following
particular case. \\

{\bf Example 2} \\

In the case 1) of Example 1, let us consider $( \Lambda, \leq ) = \mathbb{N}$, and take the
following sequence of finite partitions of $\mathbb{N}$ \\

(2.21)~~~ $ \mathbb{N} \ni \lambda \longmapsto \Delta_\lambda \in {\cal FP} ( \mathbb{N},
                                                              {\cal P} ( \mathbb{N} ) )$ \\

where \\

(2.22)~~~ $ \begin{array}{l}
               \Delta_0 = \{~ \mathbb{N} ~\} \\
               \Delta_1 = \{~ \{ 0 \}, \{ 1, 2, 3, \ldots \} ~\} \\
               \Delta_2 = \{~ \{ 0 \}, \{ 1 \}, \{ 2, 3, 4, \ldots \} ~\} \\
               \Delta_3 = \{~ \{ 0 \}, \{ 1 \}, \{ 2 \}, \{ 3, 4, 5, \ldots \} ~\} \\
               \ldots \ldots \ldots \ldots \ldots \ldots
            \end{array} $ \\

thus clearly $\Delta_0 \leq \Delta_1 \leq \Delta_2 \leq \Delta_3 \leq \ldots$ \\

Now if we take $x = 0 \in \mathbb{N} = X$, then, see (2.12) \\

(2.23)~~~ $ \Delta_\lambda ( x ) = \{~ A_{\lambda,x} =
     \{ \lambda, \lambda + 1, \lambda +2, \ldots \} ~\},~~~ \lambda \in \Lambda $ \\

Therefore (2.13) {\it fails} to hold, since obviously \\

(2.24)~~~ $ \bigcap_{\lambda \in \Lambda} A_{\lambda,x} = \phi $ \\

On the other hand, regarding (2.14), in view of (2.20), (2.23), as well as (2.16) - (2.19), we
obtain \\

(2.25)~~~ $ \underleftarrow{\lim}_{\, \lambda \,\in\, \Lambda}~ \Delta_\lambda(x) =
                        \{~ ( A_{\lambda,x} ~|~ \lambda \in \Lambda ) ~\} \neq \phi $ \\

In the case 2) of Example 1, for $x \in \mathbb{N}$, we have, see (1.6) \\

(2.26)~~~ $ \Delta_\lambda ( x ) = \{~ A_{\lambda,x} ~\} $ \\

where \\

(2.27)~~~ $ A_{\lambda,x} ~=~ \begin{array}{|l}
                                 \{ x \} ~~~\mbox{if}~~ x < \lambda \\ \\
                                 \{ \lambda, \lambda + 1, \lambda +2, \ldots \} ~~~\mbox{if}~~ x \geq \lambda
                         \end{array} $ \\

thus (2.13) will this time hold, since (2.26), (2.27) obviously yield for $x \in \mathbb{N}$ \\

(2.28)~~~ $ \bigcap_{\lambda \in \Lambda} A_{\lambda,x} = \{ x \} \neq \phi $ \\

As for (2.14), the relations (2.26), (2.27) applied to (2.20) give \\

(2.29)~~~ $ \underleftarrow{\lim}_{\, \lambda \,\in\, \Lambda}~ \Delta_\lambda(x) =
                        \{~ ( A_{\lambda,x} ~|~ \lambda \in \Lambda ) ~\} \neq \phi $ \\

which in view of (2.26) - (2.28) means essentially that \\

(2.30)~~~ $ \underleftarrow{\lim}_{\, \lambda \,\in\, \Lambda}~ \Delta_\lambda(x) =
                                                           \{~ \{ x \} ~\} \neq \phi $ \\

Lastly, in the case 3) of Example 1, we have, see (1.7) \\

(2.31)~~~ $ \Delta_\lambda ( x ) = \{~ A_{\lambda,x} ~\},~~~ x \in \mathbb{N} $ \\

where \\

(2.32)~~~ $ A_{\lambda,x} ~=~ \begin{array}{|l}
                                 \{ x_* \} ~~~\mbox{if}~~ x_* < \lambda \\ \\
                                 \{ \lambda, \lambda + 1, \lambda +2, \ldots \}
                                 ~~~\mbox{if}~~ x_* \geq \lambda
                              \end{array} $ \\

hence (2.13) holds again, since \\

(2.33)~~~ $ \bigcap_{\lambda \in \Lambda} A_{\lambda,x} = \{ x_* \} \neq \phi $ \\

while (2.14) takes the form \\

(2.34)~~~ $ \underleftarrow{\lim}_{\, \lambda \,\in\, \Lambda}~ \Delta_\lambda(x) =
                        \{~ ( A_{\lambda,x} ~|~ \lambda \in \Lambda ) ~\} \neq \phi $ \\

which in view of (2.32) means essentially that \\

(2.35)~~~ $ \underleftarrow{\lim}_{\, \lambda \,\in\, \Lambda}~ \Delta_\lambda(x) =
                                                       \{~ \{ x_* \} ~\} \neq \phi $ \\

{\bf Remark 2} \\

The three instances in Example 2 above, with their respective versions (2.25), (2.29),
(2.30), (2.34) and (2.35) of problem (2.14) as formulated in Problem 2, can give a motivation
for the use of the {\it inverse limits} in Ergodic Theory. Indeed, in each of these three
cases, the corresponding inverse limits reflect in a nontrivial manner obvious ergodic
properties of the specific mappings $T$ involved. \\
In this regard, the relevance of the inverse limit is particularly clear in the first instance
in Example 2, namely, when $T : \mathbb{N} \longrightarrow \mathbb{N}$ is the usual shift, and
when problem (2.13), as formulated in Problem 2, has a solution in (2.24) which does not give
much information about $T$, since the same relation may be obtained for many other mappings of
$\mathbb{N}$ into itself. \\
On the other hand, the inverse limit in (2.25) does give an information which is clearly
more revealing about the specific feature of $T$. \\

Of course, in analyzing the ergodic features of mappings $T$ of $\mathbb{N}$ into itself, one
can use a variety of other {\it refinement chains}, than the particular one in (2.21),
(2.22). We consider next such an example of a different refinement chain in the case 1) of
Example 1. \\

{\bf Example 3} \\

Let $X = \mathbb{N},~ \Sigma = {\cal P} ( \mathbb{N} )$ and consider the mapping $T :
\mathbb{N} \longrightarrow \mathbb{N}$ given by the usual shift $T ( x ) = x + 1,~ x \in
\mathbb{N}$. \\

Let ${\cal U}$ be a {\it free ultrafilter} on $X = \mathbb{N}$ which, we recall, means a
filter with the following two properties \\

(2.36)~~~ $ \forall~~ A \subseteq \mathbb{N} ~:~
           \mbox{either}~ A \in {\cal U}, ~\mbox{or}~ \mathbb{N} \setminus A \in {\cal U} $ \\

(2.37)~~~ $ \bigcap_{\,U \,\in\, {\cal U}}\, U = \phi $ \\

These two conditions imply that \\

(2.38)~~~ $ \forall~~ U \in {\cal U} ~:~ U ~\mbox{is infinite} $ \\

Furthermore, we also have that \\

(2.39)~~~ $ \exists~~ U \in {\cal U} ~:~ \mathbb{N} \setminus U ~\mbox{is infinite} $ \\

For $U \in {\cal U}$, let us consider the set of finite partitions $\Delta$ of $\mathbb{N}$
which contain $U$, that is, given by \\

(2.40)~~~ $ {\cal FP}_U ( \mathbb{N}, {\cal P} ( \mathbb{N} ) ) =
                 \{~ \Delta \in {\cal FP} ( \mathbb{N}, {\cal P} ( \mathbb{N} ) ) ~~|~~
                                                  U \in \Delta ~\} $ \\

Further, let us consider \\

(2.41)~~~ $ \begin{array}{l}
             {\cal FP}_{\cal U}~ ( \mathbb{N}, {\cal P} ( \mathbb{N} ) ) =
              \bigcup_{\,U \,\in\, {\cal U}}~ {\cal FP}_U
                       ( \mathbb{N}, {\cal P} ( \mathbb{N} ) ) = \\ \\
              ~~~~~~= \{~ \Delta \in {\cal FP} ( \mathbb{N}, {\cal P} ( \mathbb{N} ) )
                         ~~|~~ \exists~~ U \in {\cal U} ~:~ U \in \Delta ~\}
             \end{array} $ \\

We shall take now $( \Lambda, \leq ) = {\cal FP}_{\cal U} ( \mathbb{N}, {\cal P} (
\mathbb{N} ) )$ endowed with the partial order $\leq$ in (2.1) which corresponds to the usual
refinement of partitions. Finally, the mapping (2.8), (2.9) will simply be the identity
mapping \\

(2.42)~~~ $ \Lambda \ni \lambda = \Delta \longmapsto \Delta \in
                   {\cal FP}_{\cal U}~ ( \mathbb{N}, {\cal P} ( \mathbb{N} ) ) $ \\

Given $U \in {\cal U},~ \Delta \in {\cal FP}_{\cal U}~ ( \mathbb{N}, {\cal P} ( \mathbb{N}
) )$, with $U \in \Delta$, as well as $x \in \mathbb{N}$, it follows easily that \\

(2.43)~~~ $ U \in \Delta ( x ) $ \\

and in fact, we have the stronger property, similar with (1.5), namely \\

(2.44)~~~ $ \Delta ( x ) = \{~ A \in \Delta ~|~ A ~\mbox{is infinite} ~\} $ \\

Now it is easy to see that, in view of (2.20), we obtain \\

(2.45)~~~ $ ( A_\lambda ~|~ \lambda \in \Lambda ) \in
          \underleftarrow{\lim}_{\, \lambda \,\in\, \Lambda}~ \Delta_\lambda(x) $ \\

where for $\lambda = \Delta \in {\cal FP}_U ( \mathbb{N}, {\cal P} ( \mathbb{N} ) )$, we
have \\

(2.46)~~~ $ \Delta_\lambda = \Delta,~~~ A_\lambda = U $ \\

hence \\

(2.47)~~~ $ \underleftarrow{\lim}_{\, \lambda \,\in\, \Lambda}~
                                       \Delta_\lambda(x) \neq \phi $ \\ \\

{\bf 3. A General Inverse Limit Ergodic Result} \\

As seen in the theorem next, the result in (2.25) is in fact a particular case of a rather
general one. \\

{\bf Theorem 1} \\

Let $( X, \Sigma, T )$ be as at the beginning of section 1. Further, let $( \Lambda, \leq ) $
be a {\it directed} partial order, together with a mapping, see (2.8) \\

$~~~~~~ \Lambda \ni \lambda \longmapsto \Delta_\lambda \in {\cal FP} ( X, \Sigma ) $ \\

which satisfies (2.9), as well as the following condition : \\

(3.1)~~~ $ \exists~~ \Lambda_0 \subseteq \Lambda ~:~ \Lambda_0 ~\mbox{ is {\it countable}
                                                          and {\it cofinal} in}~ \Lambda $ \\

Then for every $x \in X$, we have \\

(3.2)~~~ $ \underleftarrow{\lim}_{\, \lambda \,\in\, \Lambda}~
                                          \Delta_\lambda(x) \neq \phi $ \\

{\bf Proof.} \\

It follows from Proposition 5 in [1, p. 198], whose conditions are satisfied, as shown
next. \\
Indeed, given $x \in X$, in view of (2.10), we have \\

(3.3)~~~ $ \Delta_\lambda(x) \neq \phi,~~~ \lambda \in \Lambda $ \\

Further, (2.16) gives for $\lambda, \lambda\,' \in \Lambda,~ \lambda \leq \lambda\,'$ the
mapping \\

(3.4)~~~ $\psi_{\lambda\,',\lambda,x} : \Delta_{\lambda\,'} \longrightarrow
                                                             \Delta_\lambda(x) $ \\

and obviously, see (2.17) - (2.19), (2.3), (2.4) \\

(3.5)~~~ $ \psi_{\lambda,\lambda,x} = id_{\Delta_\lambda(x)} $ \\

while for $\lambda, \lambda\,', \lambda\,'' \in \Lambda,~ \lambda \leq \lambda\,' \leq
\lambda\,''$ we have \\

(3.6)~~~ $ \psi_{\lambda\,',\lambda,x} \circ \psi_{\lambda\,'',\lambda\,',x} =
                                                               \psi_{\lambda\,'',\lambda,x} $ \\

Lastly, the mappings (3.3) are surjective. Indeed, let $A \in \Delta_\lambda(x)$, then we have
to find $A\,' \in \Delta_{\lambda\,'}(x)$, such that \\

(3.7)~~~ $ \psi_{\lambda\,',\lambda,x} ( A\,' ) = A $ \\

But (1.3) yields \\

(3.8)~~~~ $ \{~ n \in \mathbb{N}_+ ~|~ T^n ( x ) \in A ~\} ~~\mbox{is infinite} $ \\

Therefore (2.1) - (2.4) will give $A\,' \in \Delta_{\lambda\,'}(x)$, such that $A\,'
\subseteq A$, which means precisely (3.7).

\hfill $\Box$ \\

{\bf Remark 3} \\

1) An important fact in Theorem 1 above is that there are {\it no} conditions whatsoever
required on the mappings $T : X \longrightarrow X$. \\

2) In general, when $\Sigma$ is uncountable - a case which is often of interest in
applications - the set ${\cal FP} ( X, \Sigma )$ of all finite partitions of $X$ with subsets
in $\Sigma$, see (1.1), will also be uncountable. Furthermore, when considered with the
natural partial order (2.1), (2.2), the set ${\cal FP} ( X, \Sigma )$ does {\it not} have a
countable cofinal subset. Therefore, in such a case one {\it cannot} take in Theorem 1 \\

(3.9)~~~ $ ( \Lambda, \leq ) = ( {\cal FP} ( X, \Sigma ), \leq ) $ \\

as the directed partial order, and instead, one has to limit oneself to smaller directed
partial orders $( \Lambda, \leq )$, namely, to those which satisfy condition (3.1). \\

3) The set $\Sigma$ can be uncountable even when $X$ is countable, since one can take, for
instance, $\Sigma = {\cal P} ( X )$, that is, the set of all subsets of $X$. \\

\end{document}